\documentstyle[12pt]{article}
\renewcommand {\c}  {\'{c}}
\newcommand {\cc} {\v{c}}
\newcommand   {\s}  {\v{s}}

\begin{document}

%
%
\thispagestyle{empty}
\begin{titlepage}
\baselineskip=24pt
\title{ INTERPOLATION BETWEEN PARA-BOSE AND PARA-FERMI STATISTICS}
\author{ Stjepan Meljanac$^{1,*}$ ,Marijan Milekovi\c $^{2,**}$ and 
Ante Perica$^{1,***}$  }
\maketitle
\bigskip
\begin{center}
{\it $^{1}$ Rudjer Bo\s kovi\c \ Institute\\Bijeni\cc ka c.54, 41001 Zagreb,
 Croatia}\\
{\it $^{2}$ Prirodoslovno-Matemati\cc ki Fakultet,Zavod za teorijsku fiziku,\\
Bijeni\cc ka c.32, 41000 Zagreb, Croatia}\\
$^{*}$ E-mail: meljanac@thphys.irb.hr\\
FAX:++385-1-428541\\
$^{**}$ E-mail: marijan@phy.hr\\
$^{***}$E-mail: perica@thphys.irb.hr
\bigskip

{\bf Abstract}
\end{center}

Using deformed Green's oscillators and Green's Ansatz,we construct a 
multiparameter interpolation between para-Bose and para-Fermi statistics of a 
given order.When the interpolating parameters $q_{ij}$ satisfy 
$|q_{ij}|<1 \, ( |q_{ij}|= 1)$,the interpolation statistics is "infinite quon"-like 
( anyon-like ).The  proposed interpolation  does not contain states 
of negative norms.

\end{titlepage}

\newpage
\section{Introduction}
\setcounter{equation}{0}
\baselineskip=24pt
The properties of presently known particles are usually described in terms of an 
effective field theory in which the field operators at space-like points commute 
or anti-commute.The choice of commutation relations for bosons and 
anti-commutation relations for fermions are uniquely determined by 
the spin-statistics theorem (Pauli principle).During the past years much effort 
has been made to construct consistent generalizations of  
 Bose and Fermi statistics.One motivation comes from  theoretical and 
experimental search for possible violation of the Pauli exclusion principle 
 in 3+1 dimensions [1].The other motivation comes from the study of 
 some phenomena in condensed matter whose dynamics is essentially 
 two-dimensional [2].\\
 Theoretically,there are two types of generalized statistics which follow  
 from the relations between  creation and annihilation operators.\\ 
 (i) $\;$One type of interpolation between 
Bose and Fermi statistics is described by a ${\it continuous}$ parameter.
The characteristic example is  
infinite (quon) statistics [3,4] and the other example is anyonic (braid) statistics 
in $2+1$ dimensions [5].

The interpolating infinite (quon) statistics [3] is 
described by the following commutation relations
\begin{equation}
a_i a_j^{\dagger} -q\, a_j^{\dagger} a_i = \delta_{ij} , 
\end{equation}
$$
\forall i,j \in I , q\in {\bf R}, \quad -1 \leq q \leq +1 ,
$$
and with the vacuum condition $a_i|0\rangle  = 0$, $ \forall i \in I$.
It has been shown that the corresponding Fock space is positive definite for 
$-1 \leq q \leq +1$ [4] and that the different permutations of a given 
multiparticle state are linearly independent, i.e. all representations of the 
symmetric group can occur .

In the lower number of dimensions ,interpolating statistics are described by 
 fractional statistics of anyons.
The anyonic-type exchange algebra of annihilation and creation anyonic operators 
$a_i$,$a^{\dagger}_i$,($i \in I$) is characterized by a continuous statistical 
parameter $\lambda$ , $\lambda \in [0,1]$ [5]:

\begin{equation}
\begin{array}{c} 

a_ia^{\dagger}_j - e^{\imath \lambda \pi sgn (i-j)}\,a^{\dagger}_ja_i 
=0 \qquad i\neq j,\qquad i,j \in I ,\\[6mm]

a_ia^{\dagger}_i - \cos(\lambda\pi) \, a^{\dagger}_ia_i = 1,\\[6mm]

a_ia_j- e^{-\imath \lambda \pi sgn (i-j)}\,a_ja_i = 0 ,\\[6mm]

\end{array}
\end{equation}
The vacuum condition is $a_i|0\rangle  = 0 $, $ \forall i \in I $.
The set $\{I\}$ can be position or momentum space,discrete or continuous.
The anyonic algebra can be obtained from the Bose algebra by mapping defined in [6].
It is now widely accepted that particles with fractional statistics are likely to play 
a role in fractional Hall effects and some other phenomena 
in condensed matter [2].

Both interpolations,anyonic and quonic,can be unified in a generalized quon 
algebra [7]
\begin{equation}
a_i a_j^{\dagger} -q_{ij} a_j^{\dagger} a_i = \delta_{ij} ,\quad q_{ij}^{\ast } 
= q_{ji}, \qquad \forall i,j \in I
\end{equation}
with the vacuum condition $a_i|0\rangle  = 0$, $ \forall i \in I$.\\
(ii)$\;$ Another type of generalized statistics is parastatistics (
para-Bose and para-Fermi statistics),proposed by Green [8].These 
were the first consistent generalizations of Bose and Fermi statistics.Parastatistics 
are characterized by a ${\it discrete}$ parameter $ p \in {\bf N}$ (order of 
parastatistics).They are described by trilinear relations:
\begin{equation}
\begin{array}{c}

[a_k,[a^{\dagger}_l,a_m]_{\pm}]=(\frac{2}{p})\delta_{kl}a_m ,\\[6mm]

[a_k,[a^{\dagger}_l,a^{\dagger}_m]_{\pm}]=(\frac{2}{p})(\delta_{kl}a^{\dagger}_m 
\pm 2 \delta_{km}a^{\dagger}_l),\\[6mm]

[a_k,[a_l,a_m]_{\pm}]=0 ,

\end{array}
\end{equation}
with the vacuum conditions $a_k|0\rangle  = 0$ , $a_k a^{\dagger}_l|0\rangle =
\delta_{kl}|0\rangle,$ $k,l \in I$ and $p \in {\bf N}$ .\\The sign + (-) corresponds 
to the para-Bose (para-Fermi) algebra.Note that the last two relations are not 
independent but follow from the first one.The corresponding Fock space has no 
states with negative norms.For $p=1$ the para-Bose (para-Fermi) algebra becomes 
the Bose (Fermi) algebra and for $p =\infty $ the para-Bose (para-Fermi) 
algebra becomes 
the Fermi (Bose) algebra .The possibility of describing  some physical models 
in solids in terms of quasi-particles which obey para-Fermi statistics has 
recently been discussed  
by Safonov [9].\\
There were a few attempts to perform a continuous interpolation between  different 
parastatistics [10].The proposed interpolating trilinear commutation relations are
\begin{equation}
[a_{i} a^{\dagger}_{j} + q\,a^{\dagger}_{j}a_{i},a_k]=\rho \delta_{jk} a_i ,
\end{equation}
with the vacuum conditions $a_k|0\rangle  = 0$ , $a_k a^{\dagger}_l|0\rangle =
\delta_{kl}|0\rangle,$ $k,l \in I$.For $\rho=-(\frac{2}{p}) q $ , \\$p \in {\bf N}$ and 
$q=+1 (-1)$ , one recovers the para-Bose (para-Fermi) algebra.However,it was shown 
that for the generic ${\it q}$ the corresponding Fock-like space contained states with 
negative norms [11].Hence,no small violation of parastatistics, described by 
Eq.(5),is allowed.Recently,Speicher [12] and two of us [7]  indicated the 
possibility of continuous interpolation between parastatistics without 
states with negative norms,but no explicit construction was presented.

In this Letter we describe a general construction of  a continuous 
interpolation between  the para-Bose and 
the para-Fermi algebra of a given order ${\it p}$.The construction is performed in such a 
way that the Fock space does not contain states with negative norms.Furthermore,
for $p=1$ and $p= (\infty)$ the construction reduces to the quonic interpolation (Eq.(1)) or to 
the anyonic  interpolation (Eq.(2)),depending on interpolating parameters $q_{ij}$.

\section{q-deformed Green's Ansatz}
\baselineskip=24pt
We start with the generalized quon algebra,Eq.(3),and use it to deform the 
Green's Ansatz for the parastatistics [8] in the following way (the algebraic structure 
of the Green's Ansatz and its q-deformed analogue is described in [13]):
\begin{equation}
\begin{array}{c} 
b^{\alpha}_i b^{\dagger \beta}_j -q_{\alpha \beta}\,
b^{\dagger \beta}_j b^{\alpha}_i = \delta_{\alpha \beta} \delta_{ij},
\quad i,j \in I, \; \alpha ,\beta = 1,2,..p\\[6mm]
q_{\alpha \beta}=q\,\Delta_{\alpha\beta}\equiv q\, (2 \delta_{\alpha \beta}-1), \quad q \in {\bf R},\quad 
-1 \leq q \leq +1,
\end{array}
\end{equation}
with the vacuum conditions $b^{\alpha}_i|0\rangle=0$,
$b^{\alpha}_i b^{\dagger \beta}_j|0\rangle = \delta_{\alpha \beta}\delta_{ij}
|0\rangle $.Note that for $|q|<1$ there are no commutation relations between 
$b^{\alpha}_i$,$b^{\beta}_j$ . \\
If $q=+1 (-1)$, the oscillators $ b_{i}^{\alpha } $ and $ b^{\dagger  \alpha}_{i}  $
 reduce to 
the ordinary Green's oscillators,leading to the para-Bose (para-Fermi) algebra,Eq.(4).
The  main point of our construction is that the Fock space corresponding to the 
commutation relation (6) does not contain  states with negative norms if 
$|q| \leq 1$ (see [14,15]).

Let us define the operators $A_i$,$A_i^{\dagger}$,$i\in I$, in the same way as in
 the ordinary Green's Ansatz [8],namely:
\begin{equation}
A_i=\frac{1}{\sqrt p} \sum_{\alpha=1}^{p} b_i^{\alpha}\, ,\qquad  
A_i^{\dagger}=\frac{1}{\sqrt p} \sum_{\alpha=1}^{p} b_i^{\alpha \dagger}.
\end{equation}
The Fock space ${\cal F}(A_i)$ created by the $A_i^{\dagger}$ operators 
can be built.It is important to note that ${\cal F}(A_i)$ is a subspace of 
the initial Fock space $\cal F$ $(b_i^{\alpha})$,and therefore the space 
${\cal F}(A_i)$ automatically does not contain states with negative norms since 
$\cal F$ $(b_i^{\alpha})$  does not contain states with negative norms.

Let us define the matrix $ {\cal A}^{(n)}$ of inner products with the 
matrix elements:
\begin{equation}
{\cal A}^{(n)}_{i_1\cdots i_n;j_1\cdots j_n}=\langle 0|A_{i_n}\cdots A_{i_1}
A_{j_1}^{\dagger}\cdots A_{j_n}^{\dagger}|0\rangle
\end{equation}
where $\{i_1\cdots i_n$;$j_1\cdots j_n \} \in I$.
If indices $i_1\cdots i_n$ are 
mutually different,the matrix ${\cal A}^{(n)}$ is an (${\it n!\times n!}$) matrix 
whose 
diagonal elements are equal to 1.
\\
We find that an arbitrary matrix element is
\begin{equation}
{\cal A}^{(n)}_{\pi (i_1\cdots i_n);\sigma (i_1\cdots i_n)}=\frac{1}{p^n}
\sum_{\alpha_1,\cdots \alpha_n} 
\prod_{a,b} q_{\alpha_a \alpha_b}. 
\end{equation}
Here,$\pi$ and $\sigma$ are elements (permutations) of the permutation group $S_n$.
The product is taken over those pairs $a , b = 1,\cdots n$,which satisfy  $a<b$ and 
$(\sigma^{-1}\pi) (a)>(\sigma^{-1}\pi )(b)$.
From the general theorem on positivity [14,15] it follows that the general matrix 
${\cal A}^{(n)}$,$n \in {\bf N}$,is positive definite for $|q|<1$,and therefore the 
states $\pi (A_{i_1}^{\dagger}\cdots A_{i_n}^{\dagger})|0\rangle $,where 
$\pi \in S_n$ and ($i_1\cdots i_n$) are 
mutually different,are linearly independent as for quons [3,4].\\
However,the operators $A_i$,$A_i^{\dagger}$,$i\in I$, do not close the commutation 
relations between themselves,except for $q=\pm 1$, corresponding to the 
ordinary parastatistics,Eq.(4). Namely,
\begin{equation}
\begin{array}{c} 
A_i A_j^{\dagger}=\delta_{ij}-q\, A_j^{\dagger} A_i + q\,K_{ij}\\[6mm]

K_{ij}=(\frac{2}{p}) \sum_{\alpha=1}^p b_j^{\alpha \dagger}b_i^{\alpha}
\end{array}
\end{equation}
and $K_{ij}$ cannot be expanded in terms of the $A_i$,$A_i^{\dagger}$ operators i.e.
 cannot be eliminated from the commutation relations except for $q=\pm 1$.
For example,
the action of the $A_i$ operator on the state 
$A_j^{\dagger} A_k^{\dagger} A_m^{\dagger}|0\rangle $ is not 
contained in the space ${\cal F}(A)$ if $|q|<1$.If $q=\pm 1$ it is contained in 
${\cal F}(A)$. 

\section{Construction of interpolation between parastatistics}
\baselineskip=24pt
The operators $A_i$,$A_i^{\dagger}$ do not close and hence do not achieve desired 
interpolation between the para-Bose and  para-Fermi algebra.Nevertheless,the 
matrices ${\cal A}^{(n)}$,Eq.(9), represent the desired interpolation between 
the corresponding matrices for the para-Bose and para-Fermi algebras.To perform our 
construction of interpolating commutation relations between 
the para-Bose and para-Fermi algebra,Eq.(4),we look for the operators $a_i$,
$a_i^{\dagger}$,$i\in I$ with closed commutation relations of the type 
\begin{equation}
a_ia_j^{\dagger}=\Gamma_{ij}(a^{\dagger},a) \qquad i,j \in I ,
\end{equation}
with the vacuum conditions 
 $
a_i |0\rangle =0$, $ a_ia_j^{\dagger}|0\rangle=\delta_{ij} |0\rangle $, 
$i,j \in I$
and where $\Gamma_{ij}$ denotes a sum of all possible independent normally  
ordered terms.Then we build a Fock-like space ${\cal F}(a)$.The main requirement 
is that the matrix of inner products in  ${\cal F}(a)$, 
${\cal A}^{(n)}_{i_1\cdots i_n;j_1\cdots j_n}=\langle 0|a_{i_n}\cdots a_{i_1}
a_{j_1}^{\dagger}\cdots a_{j_n}^{\dagger}|0\rangle $, 
 should be identical to the 
corresponding matrix elements,Eq.(9). This requirement ensures that the interpolation 
between para-Bose and para-Fermi statistics, based on Eq.(9), is continuous and 
that the corresponding Fock space ${\cal F}(a)$ does not contain states of 
negative norms 
(since all matrices ${\cal A}^{(n)}$ are positive definite [14,15]).
This procedure is well defined [16].Namely,
Eq.(11) can be expanded as
$$
\Gamma_{ij} \equiv a_ia_j^{\dagger}=\delta_{ij} + C^{(ij)}a_j^{\dagger}a_i +
\sum_k\sum_{\pi,\sigma \in S_2}C^{(ij;k)}_{\pi,\sigma}\pi (a_j^{\dagger}a_k^{\dagger})
\sigma (a_k a_i)+ \cdots
$$
\begin{equation}
\cdots + \sum_{k_1,\cdots k_n} \sum_{\pi,\sigma \in S_{n+1}}
C^{(ij;k_1\cdots k_n)}_{\pi,\sigma}\pi (a_j^{\dagger}a_{k_n}^{\dagger}\cdots 
a_{k_1}^{\dagger}) \sigma (a_{k_1}\cdots a_{k_n}a_i)
\end{equation}
The summation  is performed over those permutations $\pi$,$\sigma \in S_{n+1}$ for which 
the states $\pi (a_{k_1}\cdots a_{k_n}a_i)$ are independent.Note that the expansion 
(12) could lead to Fock space with states with negative norms.However, the negative norm states do not appear,
 owing to our requirement.
The unknown 
coefficients $C_{\pi,\sigma}$can be uniquely determined [16].Moreover,the commutation 
relations (11) and (12) imply that 
$a_i {\cal F}(a) \subset {\cal F}(a)$ and 
vice versa, i.e.
\begin{equation}
\begin{array}{c}
a_ia_j^{\dagger}|0\rangle=\delta_{ij} |0\rangle ,\\[6mm]

a_ia_{i_1}^{\dagger}a_{i_2}^{\dagger}|0\rangle =\delta_{ii_1} 

a_{i_2}^{\dagger}|0\rangle + \Phi _{i_1i_2;i_1}^i\, \delta_{ii_2}a_{i_1}^{\dagger} |0\rangle ,\\[6mm]
 
a_{i}a_{i_1}^{\dagger} \cdots a_{i_n}^{\dagger}|0\rangle = 
\sum_{k=1}^{n} \sum_{\pi \in S_{n-1}} \Phi _{i_1 \cdots i_n;\pi (i_1 \cdots \not \! i_k 
\cdots i_n)}^{i} \, 
\pi (a_{i_1}^{\dagger}\cdots \not \! 
a_{i_k}^{\dagger}\cdots a_{i_n}^{\dagger}) |0\rangle ,\\[6mm]
\end{array}
\end{equation}
where the summation on the RHS is performed over the linearly independent states.The 
slash denotes the omission of the corresponding operator $a_{i_k}^{\dagger}$.
One easily finds the coefficients $\{ \Phi \}$ as  
\begin{equation}
\{ \Phi \} =\sum_{\sigma \in S_{n-1}} [{\cal A}^{(n-1)}_{\pi,\sigma}]^{-1}\, 
{\cal A}^{(n)}_{(k,\sigma);id} .
\end{equation}
Here, $(k,\sigma)$ denotes the indices $\{ k,\sigma(1),\cdots \sigma(k-1),\sigma(k+1),
\cdots \sigma(n) \}$.
We point out that  Eq.(12) represents recurrent relations for $C_{\pi,\sigma}$,
which always have a unique solution since the determinant of this linear system is 
always  regular.\\
The solution of recurrent relations ,Eq.(12),to the second order in the operators 
$a $ and $a^{\dagger}$ is:
\begin{equation}
a_ia_j^{\dagger}=\delta_{ij} + q(\frac{2}{p}-1)\, a_j^{\dagger}a_i + 
\frac{8p(p-1)q^3}{[p^2 - (p-2)^2 q^2]^2}
 \sum_{k\in I} [Y_{jk}]^{\dagger}[Y_{ik}] +\cdots
\end{equation}
with $Y_{ik}=a_ia_k - q(\frac{2}{p}-1)\,a_k a_i$ , $ -1 \leq q\leq 1$ and $ p \in {\bf N}$.
The limiting cases, $p=1$ and $p=\infty $ correspond to the quon interpolation between
 the Bose and 
Fermi algebra [3,4].\\
If $q=\varepsilon = \pm 1$ and for $ p \in {\bf N}$, we find
\begin{equation}
a_ia_j^{\dagger}=\delta_{ij} + \varepsilon \,(\frac{2}{p} - 1) a_j^{\dagger}a_i 
+\frac{\varepsilon p}{2(p-1)}\sum_{k\in I} [Y_{jk}]^{\dagger}[Y_{ik}]+\cdots,
\end{equation}
with $Y_{ik}=a_ia_k - \varepsilon (\frac{2}{p} - 1)\,a_k a_i$.
For $\varepsilon =+ 1 (- 1)$ this corresponds to the ordinary para-Bose 
(para-Fermi) commutation relations [4].It is easy to show that for $p=1$ the above 
equation reduces to the Bose (Fermi) algebra for $\varepsilon = +1 (-1)$,
respectively, since $Y_{ik}=a_ia_k - \varepsilon \,a_k a_i \equiv 0$ and these terms 
do not appear in the expansion (16).
The  similar argument holds for $p=\infty$.\\
Hence,we have constructed an infinite (quon) statistics interpolating between 
the para-Bose and para-Fermi algebra of the $p^{th}$ order.The norms of all Fock states are 
positive definite for $|q|<1$.Note that, even for the ordinary parastatistics,
Eq. (12) contains an infinite set of terms on the RHS [16].

\section{Multiparametric deformation of Green's Ansatz}
\baselineskip=24pt
We point out that our construction can be extended to an arbitrary 
multiparametric interpolation between Green's 
oscillators,Eq.(6), i.e. between the para-Bose and para-Fermi algebras.In this case 
we write Eq.(6)
 as
\begin{equation}
b^{\alpha}_i b^{\dagger \beta}_j -q_{i\alpha,j\beta}\, 
b^{\dagger \beta}_j b^{\alpha}_i = \delta_{\alpha \beta} \delta_{ij},
\end{equation}
$$
q_{i\alpha,j\beta}=q_{ij}\,\Delta_{\alpha \beta},
$$
where $q_{ij}^{\ast } = q_{ji}$ and $|q_{ij}|\leq 1$.\\
The generic matrices ${\cal A}^{(n)}$,Eq.(8),  can be easily 
calculated (see Eq.(9)):
\begin{equation}
{\cal A}^{(n)}_{\pi (i_1\cdots i_n);\sigma (i_1\cdots i_n)}=\frac{1}{p^n}\,
\sum_{\alpha_1,\cdots \alpha_n}\, (\prod_{a,b} q_{i_aj_b})
(\prod_{a,b} \Delta_{\alpha_a\alpha_b}),
\end{equation}
$$
\Delta_{\alpha_a\alpha_b} =2\delta_{\alpha_a\alpha_b}-1,
$$
where the products are taken over those pairs $a , b = 1,\cdots n $,which satisfy 
 $a<b$ and $(\sigma^{-1}\pi) (a)>(\sigma^{-1}\pi )(b)$.\\
For example,${\cal A}^{(2)}$ is given by
$$
{\cal A}^{(2)}= \left (
\begin{array}{cc}
1 & q_{ij}(\frac{2}{p}-1)\\
q_{ij}^{\ast }(\frac{2}{p}-1) & 1
\end{array}
\right).
$$

The multiparametric interpolation between para-Bose and para-Fermi oscillators to 
the second order in $a_i$, $a^{\dagger}_j$ is
\begin{equation}
\Gamma_{ij} \equiv a_ia^{\dagger}_j=\delta_{ij} + 
q_{ij}\,(\frac{2}{p}-1)\,a^{\dagger}_j a_i + 8p(p-1)\,\sum_{k\in I}Q^{ji;k}\,
[Y_{jk}]^{\dagger}[Y_{ik}]+ \cdots ,
\end{equation}
$$
Q^{ji;k}=\frac{q_{ji}q_{jk}q_{ik}}
{[p^2-(p-2)^2|q_{kj}|^2][p^2-(p-2)^2|q_{ki}|^2]},
$$
 $$
 Y_{ik}=a_ia_k-q_{ki}\,(\frac{2}{p}-1)a_ka_i.
 $$
${\it Remark}$\\
Speicher [12]  suggested ${\it q}$-para-Bose (${\it q}$- para-Fermi ) fields defined through 
the q-deformed Green's oscillators
\begin{equation}
\begin{array}{c} 
b^{\alpha}_i b^{\dagger \alpha}_j -\varepsilon \, 
b^{\dagger \alpha}_j b^{\alpha}_i =  \delta_{ij},\\[4mm]

b^{\alpha}_i b^{\dagger \beta}_j= \varepsilon q \, b^{\dagger \beta}_j 
b^{\alpha}_i , \qquad \alpha \neq \beta.
\end{array}
\end{equation}
For $q=+1$ and $\varepsilon=+1 (-1)$, this corresponds to Bose ( Fermi ) 
oscillators.For $q=-1$, one recovers para-Bose (para-Fermi) oscillators when 
$\varepsilon=+1 (-1)$, respectively.This is also a special case of  Eq.(17) 
 with
$$
q_{i\alpha ,j \beta}=\varepsilon [ (1-q)\delta_{\alpha \beta} +q ],
$$
and it interpolates between the ${\it p}$  Bose (Fermi) oscillators and 
${\it p}$  
para-Bose (para-
Fermi) oscillators, respectively.


\section{Anyonic deformation of Green's Ansatz}
\baselineskip=24pt

In the same way we can construct an anyonic-like interpolation between Green's 
oscillators, i.e. between the para-Bose and para-Fermi algebras of a given 
order $p \in {\bf N}$.In this case,in  Eq.(6) we replace $q_{\alpha \beta }$ 
with
\begin{equation}
q_{\alpha \beta} \rightarrow q_{i\alpha ,j \beta}= q_{ij}\,
\Delta_{\alpha \beta},
\end{equation}
$$
q_{ij}=\left \{ \begin{array}{ll}
e^{\imath \varphi_{ij}} & \mbox{if $i\neq j$}\\
\cos(\lambda\pi) & \mbox{if $i=j$}
\end{array}
\right.
$$
where $\varphi_{ij} \in {\bf R}$.\\
The generic matrix ${\cal A}^{(n)}$  can be 
obtained from Eq.(18).
The interpolating commutation relations (to the lowest order in $a_i$, 
$a^{\dagger}_j$) are given by Eq.(19) with $q_{ij}$ given in Eq.(21). 
For $p=1$ and $p=\infty$, they reduce to the anyonic algebra ,Eq.(2).
For $\varphi_{ij}= 0 \, (\pi) $, the relations (19,21) reduce to the ordinary 
para-Bose (para-Fermi) 
relations.\\
Contrary to the quonic interpolation of the preceding section,in this case 
there exists a 
well-defined mapping of the Bose algebra (of the Jordan-Wigner type,see[6]) leading to the 
algebra (17) (with $q_{\alpha\beta}$ given by Eq.(21)),that is
\begin{equation}
b_i^{\alpha} =e^{\imath \sum_j c_{ij}N_j \, +\imath \mu \pi \sum_{\beta}
\theta _{\alpha \beta } N_{\beta}}\, B_i^{\alpha}\, 
\sqrt {\frac{[N_{i\alpha}]_{\omega}}
{N_{i\alpha}}}
\end{equation}
$$
[N_{i\alpha}]_{\omega}=\frac{\omega^{N_{i\alpha}}-1}{\omega -1},
\qquad \omega =-\cos \lambda \pi \cos \mu \pi,
$$
$$
\varphi_{ij}=c_{ij} - c_{ji},
$$
where $\theta_{AB}$ is a step function ($\theta_{AB} = 1 \,(0)$ if $A>B\, 
(A \leq B)$),
$N_j=\sum_{\alpha} N_{j\alpha}$,\\$N_{\alpha}=\sum_j N_{j\alpha}$ are number 
operators and 
$B_i^{\alpha}$ are bosonic operators,which satisfy $[N_{i\alpha},B_j^{\beta}]=
-\delta_{\alpha \beta}\delta_{ij}\,B_i^{\alpha}$.
The resulting algebra is
$$
b_i^{\alpha}b_j^{\dagger \beta}-e^{\imath \varphi_{ij}}
e^{\imath \mu \pi sgn (\alpha 
- \beta)}\, b_j^{\dagger \beta}b_i^{\alpha}=0, \qquad (i,\alpha)\neq (j,\beta),
$$
$$
b_i^{\alpha}b_i^{\dagger\alpha}- \omega \,  b_i^{\dagger\alpha}b_i^{\alpha}=1,
$$
\begin{equation}
b_i^{\alpha}b_j^{\beta}=e^{-\imath \varphi_{ij}}
e^{-\imath \mu \pi sgn(\alpha-\beta)}\, b_i^{\beta}b_i^{\alpha},
\end{equation}
The algebra (17)(with $q_{\alpha\beta}$ given by Eq.(21)) is reproduced with $\mu=1$.
However, additional commutation relations between 
$b_i^{\alpha}$ and $b_i^{\beta}$ emerge in the anyonic case \\( $|q_{ij}|=1$),in contrast to 
 the case
 $|q_{ij}|<1$.Note also that there is  another anyonic interpolation between 
 Green's oscillators characterized by mapping,Eq.(22),with the choice of parameters 
 $c_{ij}=0$ and $\mu \in [0,1]$.


\section{Conclusion}
\baselineskip=24pt

We have considered multiparametric interpolations between the para-Bose and
 para-Fermi oscillator algebras of a given order $p \in {\bf N}$,Eqs.(4).
 The construction has been performed using deformed Green's oscillators and 
 requireing that the corresponding Fock space does not contain negative norm 
 states.
 When the interpolating parameters $q_{ij}$ in Eq.(17) satisfy $|q_{ij}|<1$,
 the interpolation between parastatistics goes through the generalized infinite (
 quon) statistics,and if $|q_{ij}|=1$,$\forall i,j \in I$, the interpolation 
 is of the anyonic type (i.e. the corresponding ${\cal A}^{(n)}$ matrices 
  are singular and 
 characteristic anyonic exchange factors appear).In the anyonic-type 
 interpolation, 
  oscillators can be obtained from para-Bose oscillators by a mapping.
 Let us remark that statistical properties are connected with the rank 
 of the matrices ${\cal A}^{(n)}$ and this is ${\it not}$ a continuous function of the 
 deformation parameters $q_{ij}$ [16].It would be interesting to explore whether such interpolating 
 statistics are compatible with  cosmological arguments presented in [17].
It is also an open question whether the anyonic-like parastatistics 
described in this Letter 
might have  applications in models of condensed matter physics and field theory.

\newpage
\baselineskip=24pt
{\bf REFERENCES}
\begin{description}

\item{[1]}
L.B.Okun, Usp.Fiz.Nauk. 158 (1989) 293 (in Russian) and references therein ;\\
O.W.Greenberg and R.N.Mohapatra, Phys. Rev. D 39 (1989) 2032;\\
R.N.Mohapatra, Phys.Lett.B 242 (1990) 407;\\
D.Miljani\c $\,$ et.al. ,Phys.Lett. B 252 (1990) 487;\\
E.Ramberg and G.E.Snow, Phys.Lett. B 238 (1990) 438. 
\item{[2]}
R.E.Prange and S.M.Girvin (eds), The Quantum Hall Effects (Springer,Berlin 1990);\\
F.Wilczek, Fractional statistics and anyon superconductivity (World Scientific,
Singapore,1990);\\
G.S.Canright and M.D.Johnson, J. Phys.A :Math.Gen. 27 (1994) 3579. 
\item{[3]}
O.W.Greenberg, Phys.Rev.Lett. 64 (1990) 705; Phys.Rev.D 43 (1991) 4111.
\item{[4]}
D.Fivel, Phys.Rev.Lett. 65 (1990) 3361;Erratum,ibid. 69 (1992) 2020;\\
M.Bozejko and R.Speicher, Commun. Math. Phys. 137 (1991) 519;\\
D.Zagier, Commun. Math. Phys. 147 (1992) 199.
\item{[5]}
J.M.Leinaas and J.Myrheim,  Nuovo Cim. 37 (1977) 1;\\
F.Wilczek, Phys.Rev.Lett. 48 (1982) 1144;\\
Y.S.Wu, Phys.Rev.Lett. 53  (1984) 111;\\
S.Forte, Rev.Mod.Phys. 64 (1992) 193;\\
V.Bardek,M.Dore\s i\c $\,$ and S.Meljanac, Phys. Rev. D 49 (1994) 3059 ;\\ 
V.Bardek,M.Dore\s i\c $\,$ and S.Meljanac, Int.J.Mod.Phys. A 9 (1994) 4185;\\
A.Liguori and M.Mintchev, Comm.Math.Phys. 169 (1995) 635 and references therein.
\item{[6]}
S.Meljanac,M.Milekovi\c $\,$ and A.Perica, Europhys.Lett. 28 (1994) 79;\\
M.Dore\s i\c ,S.Meljanac and M.Milekovi\c, Fizika 3 (1994) 57.
\item{[7]}
 S.Meljanac and A.Perica, J. Phys.A :Math.Gen. 27 (1994) 4737;
 \\ Mod.Phys.Letters
 A 9 (1994) 3293.
\item{[8]}
H.S.Green, Phys.Rev. 90 (1953) 170;\\
Y.Ohnuki and S.Kamefuchi, Quantum field theory and parastatistics 
( University of Tokio Press, Tokio, Springer, Berlin, 1982).
\item{[9]}
V.L.Safonov, Phys.Status Solidi B 167 (1991) 109; 
Mod.Phys.Letters B 8 (1994) 1195.
\item{[10]}
A.Yu.Ignatiev and V.A.Kuzmin, Yad.Fiz. 46 (1987) 786;\\
O.W.Greenberg and R.N.Mohapatra, Phys. Rev. Lett. 59 (1987) 2507;\\
A.B.Govorkov, Nucl.Phys. B 365 (1991) 381;Fiz.El.Chastits At.Yad. 4 (1993) 1341 
(in Russian) and references therein.
\item{[11]}
A.B.Govorkov, Phys.Lett. A 137 (1991) 381;\\
O.W.Greenberg and R.N.Mohapatra, Phys. Rev. Lett. 62 (1989) 712;
\item{[12]}
R.Speicher, Lett.Math.Phys. 27 (1993) 97.
\item{[13]}
T.D.Palev, Algebraic structure of the Green's ansatz and its q-deformed 
analogue, (INRNE-TH-94/4 preprint).
\item{[14]}
S.Meljanac and D.Svrtan, Multiparametric extension of Zagier's results on 
infinite statistics models, (IRB-TH-95 preprint,unpublished).
\item{[15]}
M.Bozejko and R.Speicher, Math.Ann. 300 (1994) 97.
\item{[16]}
S.Meljanac and M.Milekovi\c , Unified view of multimode algebras with Fock-like 
representations, (IRB-TH-95 preprint,to appear in Int.J.Mod.Phys. A).
\item{[17]}
J.W.Goodison and D.J.Toms, Phys.Rev.Lett. 71 (1993) 3240;\\
V.Bardek,S.Meljanac and A.Perica, Phys.Lett. B 338 (1994) 20.

\end{description}
\end{document}